\documentclass[reqno,10pt]{article}
\usepackage[english]{babel}
\usepackage{amsfonts}
\usepackage{amsmath}
\usepackage{graphicx}
\usepackage{amssymb}
\usepackage{amsthm}
\newcommand{\Real}{\mathbb{R}}

\newcommand{\inter}{\mathrm{int}}

\newcommand{\dom}{\mathrm{dom}}
\newcommand{\Fix}{\mathrm{Fix}}
\newcommand{\cover}[1]{\stackrel{#1}{\Longrightarrow}}
\newcommand{\backcover}[1]{\stackrel{#1}{\Longleftarrow}}
\newcommand{\bicover}[1]{\stackrel{#1}{\Longleftrightarrow}}
\newcommand{\st}{\ |\ }
\def\oarrow{{\mathrel{\longrightarrow\mkern-25mu\circ}\;\;}}
\newcommand{\comment}[1]{}
\newtheorem{thm}{Theorem}[section]

\newtheorem{lem}[thm]{Lemma}

\newtheorem{defn}[thm]{Definition}
\newtheorem{rem}[thm]{Remark}
\begin{document}
\begin{center}
{\bf \Large  Heteroclinic Connections between Periodic Orbits in
Planar Restricted Circular Three Body Problem - Part II}

\vskip 2\baselineskip {\large Daniel Wilczak}\footnote{research
supported by Polish State Committee for Scientific
Research grant 2 P03A 041 24 %  grant Mrozka
}  \\
  WSB -- NLU, Faculty of Computer Science,\\
  Department of Computational Mathematics,\\
  Zielona 27, 33-300 Nowy S\c{a}cz, Poland\\
  and \\
  Jagiellonian University, Institute of Computer Science\\
  Nawojki 11, 30-072  Krak\'ow, Poland\\
  e-mail: dwilczak@wsb-nlu.edu.pl

\vskip \baselineskip {\large and}\vskip \baselineskip

    {\large Piotr Zgliczy\'nski}  \\
    Jagiellonian University, Institute of Mathematics, \\
    Reymonta 4, 30-059 Krak\'ow, Poland\\
    e-mail: zgliczyn@im.uj.edu.pl

\vskip \baselineskip \today
\end{center}

\begin{abstract}
We present a method for proving the existence of symmetric
periodic, heteroclinic or homoclinic orbits in dynamical systems
with the reversing symmetry. As an application we show that the
Planar Restricted Circular Three Body Problem (PCR3BP)
corresponding to the Sun-Jupiter-Oterma system possesses an
infinite number of symmetric periodic orbits and homoclinic orbits
to the Lyapunov orbits. Moreover, we show the existence of
symbolic dynamics on six symbols for PCR3BP and the possibility of
resonance transitions  of the comet. This extends earlier results
by Wilczak and Zgliczynski\cite{WZ}.
\end{abstract}
%
%------------------------------------------------------------------------------------------------
%
\section{Introduction.}

 The Planar Restricted
Circular Three Body Problem (PCR3BP) has attracted much attention
of scientists. In particular, some transport properties of this
system may be applied in space mission design (see \cite{KLMR} and
references given there). The problem has been studied by Koon, Lo,
Marsden and Roos in \cite{KLMR}, where  the numerical evidence of
the resonance transitions for the PCR3BP for the parameter values
corresponding to the Sun-Jupiter-Oterma system is presented. The
rigorous proof of some facts discovered in \cite{KLMR} was given
by Wilczak and Zgliczynski in \cite{WZ} (see also Stoffer and
Kirchgraber paper \cite{SK}).

  In the present paper, as in \cite{WZ},
 we restrict our attention to the following parameter values for PCR3BP $C=3.03$,
 $\mu=0.0009537$ - the parameter values for {\em Oterma} comet in the Sun-Jupiter system
(see \cite{KLMR}).  For this parameter values there are two
hyperbolic periodic orbits $L_1^*$ and $L^*_2$, the Liapunov
orbits, around the libration points $L_1$ and $L_2$, respectively.
In \cite{WZ} and \cite{SK} it was proven that  there exist
homoclinic solutions to both $L_1^*$ and $L_2^*$ periodic orbits
and a pair of heteroclinic connections between them in both
directions.

In this paper we present the proof of the following facts:
\begin{itemize}
\item The PCR3BP possesses two pairs of homoclinic orbits both to
 $L_1^*$ and $L_2^*$. These homoclinic orbits are
geometrically different. Informally speaking they are close to
different resonances, namely 3:2, 5:3 for the orbits homoclinic to
$L_1^*$ and 1:2, 2:3 for the orbits homoclinic to $L_2^*$.
Moreover, it is possible for a comet to move between these four
resonances in arbitrary order.

\item The PCR3BP possesses an infinite number  of geometrically
different symmetric periodic and homoclinic orbits.
\end{itemize}

Let describe now what constitute a new element in the present
paper. While the numerical evidence of three of the above
mentioned homoclinic orbits is given in \cite{KLMR}, the 2:3
homoclinic orbit appear to be a new one. The technique of the
proof of the existence of these homoclinic orbits and the symbolic
dynamics  is the same as in \cite{WZ}, i.e. combines topological
tools (covering relations) with rigorous numerics.

The main novelty of the present paper, when compared to
\cite{KLMR,WZ}, is the existence of an infinite number  of
geometrically different symmetric periodic and homoclinic orbits.
While the numerical evidence of the simplest symmetric homoclinic
orbits is given in \cite{KLMR} the numerical method used there
cannot yield the existence of an infinite number of them even with
the help of validated numerics. In this paper we use the
topological method introduced recently by Wilczak in \cite{W2} and
developed later by both authors in \cite{WZ1}. For the purpose of
this introduction we briefly describe the main points of  method
for symmetric periodic points. The method is based on two
observations:
\begin{itemize}
\item to detect symmetric periodic orbits for map $P$ (a Poincar\'e map) with
an reversing symmetry $R$ (in the PCR3BP case composition of a
suitable reflection and the time inversion) it is enough to look
for intersections of $\Fix(R)=\{x \: | \: R(x)=x \}$ with
$P^k(\Fix(R))$. This is the Fixed Set Iteration method
\cite{la,la1} (also known as DeVogelaere method \cite{dv}). Any
point from such  intersection  give rise to $2k$-periodic point.
\item covering relations give some  control  of  pieces of
$P^k(\Fix(R))$, which make it possible to prove that $P^k(\Fix(S))
\cap \Fix(R)$ is nonempty for $k$ sufficiently large and that the
period of the periodic point is indeed equal to $2k$.
\end{itemize}

The paper is organized as follows. In Section~\ref{sec:descr} we
recall the PCR3BP and its properties. In
Section~\ref{sec:homoclinic} the proof of the existence of a new
pair of homoclinic orbits both in exterior and interior regions to
the Lyapunov orbits $L_1^*$, $L_2^*$ is presented. In
Section~\ref{sec:symdyn} the symbolic dynamics on six symbols is
established. We also discuss the resonance transitions there. In
Section~\ref{sec:symmetry} the existence of symmetric periodic and
homoclinic orbits is proven.

Throughout the paper we will use the definitions and notations
from \cite{WZ}.

\section{Short description of the system.}\label{sec:descr}

We follow papers \cite{KLMR,WZ} and  use the notation introduced
there.

Let $S$ and $J$ be two bodies called Sun and Jupiter, of masses
$m_s=1-\mu$ and $m_j=\mu$, $\mu \in (0,1)$, respectively. They
rotate in the plane in circles counter clockwise about their
common center and with angular velocity normalized as one. Choose
a rotating coordinate system,  so that origin is at the center of
mass and the Sun and Jupiter are fixed on the $x$-axis at
$(-\mu,0)$ and $(1-\mu,0)$ respectively. In this coordinate frame
the equations of motion of a massless particle called the comet or
the spacecraft under the gravitational action of Sun and Jupiter
are (see \cite{KLMR} and references given there)
\begin{equation}
  \ddot{x} - 2\dot{y}=\Omega_x(x,y), \qquad
  \ddot{y} + 2\dot{x}=\Omega_y(x,y),  \label{eq:PCR3BP}
\end{equation}
where
\begin{eqnarray*}
  \Omega(x,y)=\frac{x^2 + y^2}{2} + \frac{1 - \mu}{r_1} + \frac{\mu}{r_2}
     + \frac{\mu(1-\mu)}{2} \\
   r_1=\sqrt{(x+\mu)^2 + y^2}  , \qquad r_2=\sqrt{(x-1+\mu)^2 + y^2}
\end{eqnarray*}
Equations (\ref{eq:PCR3BP}) are called the equations of the planar
circular restricted three-body problem (PCR3BP). They have a first
integral called the {\em Jacobi integral}, which is given by
\begin{equation}
  C(x,y,\dot{x},\dot{y})= - (\dot{x}^2 + \dot{y}^2) +
  2\Omega(x,y).
\end{equation}

We  consider PCR3BP on the hypersurface
\begin{equation*}
  \mathcal{M}(\mu,C)=\{ (x,y,\dot{x},\dot{y}) \: | \: C(x,y,\dot{x},\dot{y})=C
  \},
\end{equation*}
and we restrict our attention to the following parameter values
$C=3.03$, $\mu=0.0009537$ - the parameter values for {\em Oterma}
comet in the Sun-Jupiter system (see \cite{KLMR}).

The projection of $\mathcal{M}(\mu,C)$ onto position space is
called a Hill's region and gives the region in the $(x,y)$-plane,
where the comet is free to move. The Hill's region for the
parameter considered in this paper is shown on Figure
\ref{fig:homoc1} in white, the forbidden region is dark. The
Hill's region consists of three regions: an interior (Sun) region,
an exterior region and Jupiter region.

As was mentioned in the Introduction we restrict our attention to
the following parameter values $C=3.03$, $\mu=0.0009537$ - the
parameter values for {\em Oterma} comet in the Sun-Jupiter system
(see \cite{KLMR}). Since we work with fixed parameter values we
usually drop the dependence of various objects defined throughout
the paper on $\mu$ and $C$, so for example
$\mathcal{M}=\mathcal{M}(\mu,C)$.

\subsection{Poincar\'e maps.}
We consider  Poincar\'e sections: $\Theta =
\{(x,y,\dot{x},\dot{y})\in \mathcal{M} \ | \ y=0\}$, $\Theta_+ =
\Theta\cap\{\dot{y}>0\}$, $\Theta_- = \Theta\cap\{\dot{y}<0\}$.

On $\Theta_\pm$ we can express $\dot{y}$ in terms of $x$ and
$\dot{x}$ as follows
\begin{equation*}
\dot{y} = \pm\sqrt{2\Omega(x,0) - \dot{x}^2-C}
\end{equation*}
Hence the sections $\Theta_\pm$ can be parameterized by two
coordinates $(x,\dot{x})$ and we will use this identification
throughout the paper. More formally, we have the transformation
$T_{\pm}:\mathbb{R}^2 \to \Theta_{\pm}$
 given by the following formula
\begin{equation*}
 T_\pm(x,\dot{x})=(x,0,\dot{x},\pm\sqrt{2\Omega(x,0) - \dot{x}^2-C}\mbox{ })
\end{equation*}
The domain of $T_\pm$ is given by an inequality $2\Omega(x,0) -
\dot{x}^2-C \geq 0$.

Let $\pi_{\dot{x}}:\Theta_\pm \longrightarrow\mathbb{R}$ and
$\pi_{x}:\Theta_\pm \longrightarrow\mathbb{R}$ denote the
projection onto $\dot{x}$ and $x$ coordinate, respectively. We
have $\pi_{\dot{x}}(x_0,\dot{x}_0)=\dot{x}_0$ and
$\pi_{x}(x_0,\dot{x}_0)=x_0$.

We will say that $(x,\dot{x}) \in \Theta_\pm$ meaning that
$(x,\dot{x})$ represents two-dimensional coordinates of a point on
$\Theta_\pm$. Analogously we  give a meaning to the statement $M
\subset \Theta_\pm$ for a set $M \subset \mathbb{R}^2$.

We define the following Poincar\'e maps between sections
\begin{eqnarray*}
  P_+: \Theta_+ \to \Theta_+ \\
  P_-: \Theta_- \to \Theta_- \\
  P_{\frac{1}{2},+}: \Theta_+ \to \Theta_- \\
  P_{\frac{1}{2},-}: \Theta_- \to \Theta_+.
\end{eqnarray*}
As a rule the sign $+$ or $-$ tells that the domain of the maps
$P_\pm$ or $P_{\frac{1}{2},\pm}$ is contained in $\Theta_\pm$ (the
same sign). Observe that
\begin{displaymath}
   P_+(x)=P_{\frac{1}{2},-} \circ P_{\frac{1}{2},+}(x), \qquad
   P_-(x)=P_{\frac{1}{2},+} \circ P_{\frac{1}{2},-}(x)
\end{displaymath}
whenever $P_+(x)$ and $P_-(x)$ are defined. These identities
express the following simple fact: to return to $\Theta_+$ we need
to cross $\Theta$ with negative $\dot{y}$ (this is
$P_{\frac{1}{2},+}$ first and then we return to $\Theta$ with
$\dot{y} >0$ (this is $P_{\frac{1}{2},-}$).

Sometimes we will drop signs in $P_\pm$ and $P_{\frac{1}{2},\pm}$,
hence $P(z)=P_{+}(z)$ if $z \in \Theta_+$ and $P(z)=P_{-}(z)$ if
$z \in \Theta_-$, a similar convention will be applied to
$P_{\frac{1}{2}}$.

\subsection{Symmetry properties of PCR3BP}
\label{subsec:symPCR3BP}
 Notice that PCR3BP has the following symmetry
\begin{equation*}
R(x,y,\dot{x},\dot{y},t)=(x,-y,-\dot{x},\dot{y},-t),
\end{equation*}
 which
expresses the following fact, if $(x(t),y(t))$ is a trajectory for
PCR3BP, then $(x(-t),-y(-t))$ is also a trajectory for PCR3BP.
From this it follows immediately that
\begin{equation}\label{eq:sym_Pf}
    \begin{split}
  \mbox{if} \quad P_\pm(x_0,\dot{x}_0)=(x_1,\dot{x}_1) \qquad \mbox{then} \quad
    P_\pm(x_1,-\dot{x}_1)=(x_0,-\dot{x}_0),  \\
  \mbox{if} \quad P_{\frac{1}{2},\pm}(x_0,\dot{x}_0)=(x_1,\dot{x_1}) \qquad \mbox{then} \quad
    P_{\frac{1}{2},\mp}(x_1,-\dot{x}_1)=(x_0,-\dot{x}_0).
    \end{split}
\end{equation}

We will denote also by $R$ the map $R:\Theta_\pm \to \Theta_\pm$
$R(x,\dot{x})=(x,-\dot{x})$ for $(x,\dot{x}) \in \Theta_\pm$. Now
Eq. (\ref{eq:sym_Pf}) can be written as
\begin{equation*}
\begin{split}
 \mbox{if} \quad P_\pm(x_0)=x_1 \qquad \mbox{then} \quad P_\pm(R(x_1))=R(x_0), \\
  \mbox{if} \quad P_{\frac{1}{2},\pm}(x_0)=x_1 \qquad \mbox{then} \quad
  P_{\frac{1}{2},\mp}(R(x_1))=R(x_0).
    \end{split}
\end{equation*}

\section{The existence of new homoclinic
orbits.}\label{sec:homoclinic}

The goal of this section is to present the proof of the existence
of new homoclinic orbits with  different resonances.

\textbf{The notion of the resonance.} We rewrite here an informal
definition of the resonance from \cite[Sec. 5.1]{KLMR}. Recall
that the PCR3BP is a perturbation of the two-body problem. Hence,
outside a small neighborhood of Jupiter, the trajectory of a comet
follows essentially a two-body orbit around the Sun. In the
heliocentric inertial frame, the orbit is nearly elliptical. The
mean motion resonance of the comet with respect to Jupiter is
equal to $a^{-3/2}$ where $a$ is the semi-major axis of this
elliptical orbit. Recall that the Sun-Jupiter distance is
normalized to be $1$ in the PCR3BP. The comet is said to be in
$p:q$ resonance with Jupiter if $a^{-3/2}\approx p/q$, where $p$
and $q$ are small integers. In heliocentric inertial frame, the
comet makes roughly $p$ revolutions around the Sun in $q$ Jupiter
periods. Observe that this definition of the resonance make also
sense for the orbits, which are non-periodic (for example orbits
homoclinic to $L^*_1$ or $L_2^*$), we just have to compute the
semi-major axis for the piece of orbit away from Jupiter. A
heuristic approach, which allows to read the resonance of an orbit
from the trajectory in the rotating frame is described in
Appendix.

\begin{figure}[htbp]
\centerline{\includegraphics[width=2.5in]{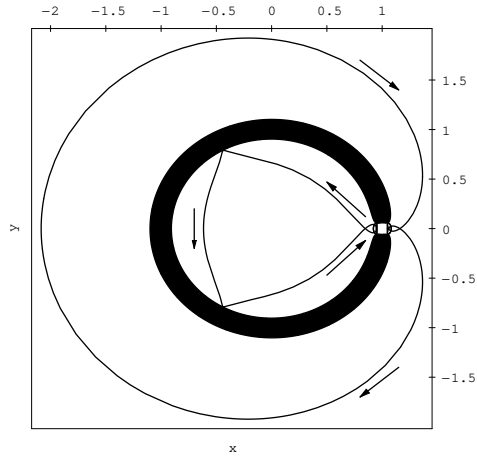}}
\caption{3:2 homoclinic orbit to $L_1^*$ Lyapunov orbit (interior
region) and 1:2 homoclinic orbit to $L_2^*$ Lyapunov orbit
(exterior region).\label{fig:homoc1}}
\end{figure}

In \cite{WZ} the following theorem was proved.
\begin{thm}{\cite[Thm.6.5,Thm.6.7]{WZ}}
\label{thm:hom} Consider  PCR3BP with $C=3.03$, $\mu=0.0009537$.
Then
\begin{itemize}
\item there exist a homoclinic orbit to the $L_1^*$ orbit (in Sun region).
This orbit is close to the $3:2$ resonance.
\item there exist a homoclinic orbit to the $L_2^*$ orbit (in exterior
region). This orbit is close to the $1:2$ resonance.
\end{itemize}
\end{thm}
These orbits are presented in Fig.~\ref{fig:homoc1}.

In this section we establish the existence of new homoclinic
connections both in exterior and interior regions. The new
homoclinic orbit in exterior region is close to the 2:3 resonance.
As was mentioned in the Introduction this orbit has been found
numerically in \cite{KLMR}, see Fig. 5.4  the and the intersection
stable and unstable manifolds of $L^*_2$ at $L=\sqrt{a}\approx
1.26$. The other new  homoclinic orbit in interior region is close
to the 5:3 resonance appears to be a new one.

\begin{figure}[htbp]
\centerline{\includegraphics[width=2.5in]{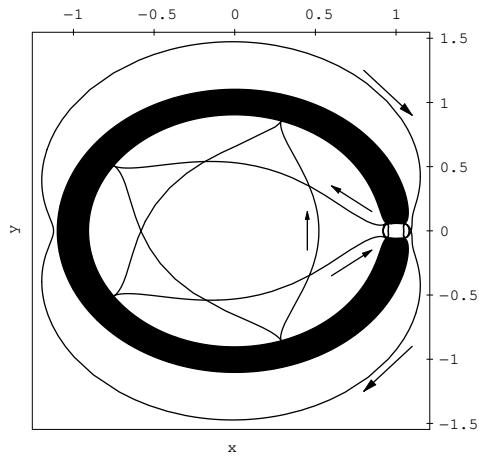}}
\caption{5:3 homoclinic orbit to $L_1^*$ Lyapunov orbit (interior
region) and 2:3 homoclinic orbit to $L_2^*$ Lyapunov orbit
(exterior region).\label{fig:homoc2}}
\end{figure}
\subsection{ The existence of the $2:3$ homoclinic orbit in the exterior region.}
We define the following h-sets $G_i=t(c_i,u_i,s_i)$, for $i=0,\ldots,4$, where
\begin{eqnarray*}
    c_0 & = & (-1.12327231155833984,0),\\
    c_1 & = & (1.093337837571255552,-0.02510094170679043584),\\
    c_2 & = & (1.047131544421841024,-0.001056187943513949696),\\
    c_3 & = & (1.08194053721089792,-2.521361165903333888\cdot10^{-5}),\\
    c_4 & = & (1.04682616720451456,-9.169345277545603072\cdot10^{-7})\\
\end{eqnarray*}
and
\begin{eqnarray*}
\begin{array}{lllclll}
 s_0 & = & (-1\cdot10^{-8},4\cdot10^{-7}), & &
 u_0 & = & -R(s_0),\\
 s_1 & = & (1\cdot10^{-7},21\cdot10^{-8}), & &
 u_1 & = & -R(s_1)/10,\\
 s_2 & = & (-1\cdot10^{-7},35\cdot10^{-8}), & &
 u_2 & = & -R(s_2)/10,\\
 s_3 & = & (-1\cdot10^{-7},23\cdot10^{-8}), & &
 u_3 & = & -R(s_3)/10,\\
 s_4 & = & (-1\cdot10^{-7},35\cdot10^{-8}), & &
 u_4 & = & -R(s_4)/4.\\
\end{array}
\end{eqnarray*}
We assume, that $G_0,G_2,G_4\subset\Theta_+$ and
$G_1,G_3\subset\Theta_-$. With a computer assistance we proved the
following
\begin{lem}
\label{lem:cov_ext_reg} The maps
\begin{eqnarray*}
  P_{\frac{1}{2},+}&:& G_0 \cup G_2 \cup G_4 \to \Theta_-, \\
  P_{\frac{1}{2},-}&:& G_1 \cup G_3 \to \Theta_+
\end{eqnarray*}
are well defined and continuous. Moreover, the following covering relations hold
\begin{eqnarray*}
G_0\cover{P_{1/2,+}}G_1\cover{P_{1/2,-}}G_2\cover{P_{1/2,+}}G_3\cover{P_{1/2,-}}G_4\cover{P_{1/2,+}}H_2^2.
\end{eqnarray*}
\end{lem}
\begin{thm}
For PCR3BP with $C=3.03$ and $\mu=0.0009537$ there exists an orbit
homoclinic to $L_2^*$  close to the $2:3$ resonance.
\end{thm}
\begin{proof}
From Lemma~\ref{lem:cov_ext_reg} and \cite[Lemma~5.6]{WZ} it
follows that
\begin{equation*}
G_0\cover{P_{1/2,+}}G_1\cover{P_{1/2,-}}G_2\cover{P_{1/2,+}}G_3\cover{P_{1/2,-}}G_4\cover{P_{1/2,+}}H_2^2\cover{P_-}H_2\cover{P_-}H_2
\end{equation*}
Note that the h-set $G_0$ is $R$-symmetric by its definition. Therefore
\begin{gather*}
    \begin{split}
    H_2=R(H_2)\backcover{P_-}R(H_2)\backcover{P_-}R(H_2^2)\backcover{P_{1/2,-}}R(G_4)\backcover{P_{1/2,+}}R(G_3)\\
    R(G_3)\backcover{P_{1/2,-}}R(G_2)\backcover{P_{1/2,+}}R(G_1)\backcover{P_{1/2,-}}R(G_0)=G_0
    \end{split}
\end{gather*}
Since $P_-$ is hyperbolic on $|H_2|$ (\cite[Lemma~5.5]{WZ}) the
assertion is a consequence of \cite[Theorem~4]{gaz}.
\end{proof}
\subsection{The existence of the $5:3$ homoclinic orbit in the interior region.}
As in the previous section we construct a chain of covering relations in order to prove the
existence of homoclinic orbit to $L_1^*$ orbit. We define h-sets
$V_i=t(c_i,u_i,s_i)$, for $i=0,\ldots,4$, where
\begin{eqnarray*}
    c_0 & = & (0.5217056203008400006,0),\\
    c_1 & = & (-0.5822638014577352639,-0.2793408708392046136),\\
    c_2 & = & (0.919204446847046941,0.004093829363524479834),\\
    c_3 & = & (0.9522506335647477061,0.0001333182992547130779),\\
    c_4 & = & (0.9208022956271231241,2.918364277340028028\cdot10^{-6})\\
\end{eqnarray*}
and
\begin{eqnarray*}
\begin{array}{lllclll}
 s_0 & = & (-1\cdot10^{-7},2\cdot10^{-7}), & & u_0 & = & -R(s_0),\\
 s_1 & = & (2\cdot10^{-8},4\cdot10^{-7}), & & u_1 & = & (3\cdot10^{-8},0),\\
 s_2 & = & (-4\cdot10^{-7},102\cdot10^{-8}), & & u_2 & = & -R(s_2)/5,\\
 s_3 & = & (-1\cdot10^{-7},365\cdot10^{-9}), & & u_3 & = & -R(s_3)/10,\\
 s_4 & = & (-1\cdot10^{-7},25733011\cdot10^{-14}), & & u_4 & = & -R(s_4)/2.\\
\end{array}
\end{eqnarray*}
We assume, that $V_0,V_2,V_4\subset\Theta_+$ and $V_1,V_3\subset\Theta_-$.
With a computer assistance we proved the following
\begin{lem}\label{lem:cov_int_reg}
The maps
\begin{eqnarray*}
  P_{\frac{1}{2},+}&:& V_0 \cup V_2 \cup V_4 \to \Theta_-, \\
  P_{\frac{1}{2},-}&:& V_1 \cup V_3 \to \Theta_+
\end{eqnarray*}
are well defined and continuous. Moreover, we have the following
chain of covering relations
\begin{eqnarray*}
V_0\cover{P_{1/2,+}}V_1\cover{P_{1/2,-}}V_2\cover{P_{1/2,+}}V_3\cover{P_{1/2,-}}V_4\cover{P_+}H_1^2.
\end{eqnarray*}
\end{lem}
\begin{thm}
For PCR3BP with $C=3.03$ and $\mu=0.0009537$ there exists an orbit
homoclinic to $L_1^*$  close to the $5:3$ resonance.
\end{thm}
\begin{proof}
Since the sets $H_1$ and $H_1^2$ are $R$-symmetric \cite[Lemma~5.6]{WZ} and \cite[Corollary~3.14]{WZ}
imply that
\begin{equation*}
H_1^2=R(H_1^2)\backcover{P_+}R(H_1)=H_1\cover{P_+}H_1.
\end{equation*}
After combining the above with Lemma~\ref{lem:cov_int_reg} we obtain
\begin{equation*}
V_0\cover{P_{1/2,+}}V_1\cover{P_{1/2,-}}V_2\cover{P_{1/2,+}}V_3\cover{P_{1/2,-}}V_4\cover{P_+}H_1^2\backcover{P_+}H_1\cover{P_+}H_1.
\end{equation*}
Note that the h-set $V_0$ is $R$-symmetric by its definition. Therefore
\begin{gather*}
    \begin{split}
    H_1\cover{P_+}H_1^2\backcover{P_+}R(V_4)\backcover{P_{1/2,+}}R(V_3)\backcover{P_{1/2,-}}R(V_2)
    \backcover{P_{1/2,+}}R(V_1)\backcover{P_{1/2,-}}R(V_0)=V_0
    \end{split}
\end{gather*}
Since $P_+$ is hyperbolic on $|H_1|$ (\cite[Lemma~5.5]{WZ}) the
assertion is a consequence of \cite[Theorem~4]{gaz}.
\end{proof}

\section{Symbolic dynamics on six symbols and resonance transitions.}\label{sec:symdyn}
\begin{figure}[htpb]
\centerline{\includegraphics[width=2.5in]{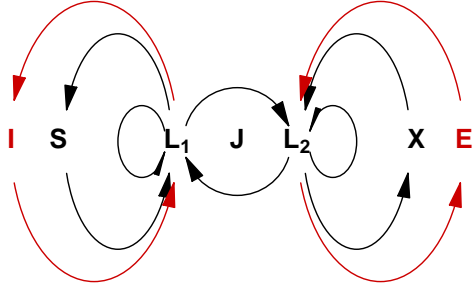}} \caption{The
graph of symbolic dynamics on six symbols.\label{fig:graf2}}
\end{figure}

As a consequence of theorems proved in \cite{WZ} and in the
previous section we obtain the existence of symbolic dynamics on
six symbols. Let $L_1$, $L_2$ denote the Lyapunov orbits regions
(see \cite{KLMR}) , $S$ and $I$ denote two parts of the Sun region
corresponding to suitable vicinities of two homoclinic orbits to
$L_1^*$. Let $X$ and $E$ denote two parts the exterior region
corresponding to suitable vicinities of two homoclinic orbits to
$L_2^*$ orbit. Schematically this situation is shown in
Fig.~\ref{fig:graf2}

In \cite{WZ} the symbolic dynamics on four symbols, i.e. $\{L_1,L_2,X,S\}$ was
established. The new homoclinic orbits allow us to include more symbols in it.
We state this result more precisely. Let
$\alpha,\beta\in\{L_1,L_2,X,E,I,S\}$ be such that there is an
arrow from $\alpha$ to $\beta$ on the graph presented in
Fig.~\ref{fig:graf2}. We define the function
\begin{equation}\label{eq:def_f}
f_{\beta,\alpha}=
\begin{cases}
    P_+, & \text{if } (\alpha,\beta)=(L_1,L_1)\\
    P_-, & \text{if } (\alpha,\beta)=(L_2,L_2)\\
    P_- \circ (P_{1/2,+} \circ  P_{1/2,-})^4 \circ P_{1/2,+} \circ  P_+, & \text{if } (\alpha,\beta)=(L_1,L_2)\\
    P_+ \circ P_{1/2,-}  \circ  (P_{1/2,+} \circ P_{1/2,-})^4 \circ P_-, & \text{if } (\alpha,\beta)=(L_2,L_1)\\
    P_+ \circ (P_{1/2,-} \circ P_{1/2,+})^2 \circ P_{1/2,-} & \text{if } (\alpha,\beta)=(S,L_1) \\
    P_+ \circ P_{1/2,-} \circ (P_{1/2,+} \circ P_{1/2,-})^2 & \text{if } (\alpha,\beta)=(L_1,S) \\
    P_-^2 \circ P_{1/2,+} \circ (P_{1/2,-} \circ P_{1/2,+})^2 &  \text{if } (\alpha,\beta)=(X,L_2)\\
    (P_{\frac{1}{2},-} \circ P_{1/2,+})^2 \circ P_{1/2,-} \circ P_{-}^2 &  \text{if } (\alpha,\beta)=(L_2,X) \\
    P_-\circ (P_{1/2,+}\circ P_{1/2,-})^2\circ P_{1/2,+}, & \text{if } (\alpha,\beta)=(E,L_2)\\
    P_{1/2,-}\circ (P_{1/2,+}\circ P_{1/2,-})^2\circ P_-, & \text{if } (\alpha,\beta)=(L_2,E)\\
    P_+^2\circ (P_{1/2,-}\circ P_{1/2,+})^2, & \text{if } (\alpha,\beta)=(I,L_1)\\
    (P_{1/2,-}\circ P_{1/2,+})^2\circ P_+^2, & \text{if } (\alpha,\beta)=(L_1,I)
\end{cases}
\end{equation}
For each symbol $\alpha\in\{L_1,L_2,X,E,I,S\}$ we define the h-set $Q_\alpha$,
where $Q_{L_1}=H_1$, $Q_{L_2}=H_2$, $Q_{S}=E_0$, $Q_X=F_0$, $Q_I=V_0$, $Q_E=G_0$.
\begin{defn}
The bi-infinite sequence $(\alpha_i)_{i\in\mathbb Z}$ is called
\emph{admissible} if for every $i\in\mathbb Z$ there is an arrow
from $\alpha_i$ to $\alpha_{i+1}$ on the graph presented in
Fig.~\ref{fig:graf2}.

The finite sequence $(\alpha_0,\alpha_1,\ldots,\alpha_n)$ is
called \emph{admissible} if for every $i=0,1,\ldots, n-1$ there is
an arrow from $\alpha_i$ to $\alpha_{i+1}$ on the graph presented
in Fig.~\ref{fig:graf2}.
\end{defn}

Let $\Gamma$ be the set of all admissible sequences
$(\alpha_i)_{i\in\mathbb Z}\in\{L_1,L_2,X,E,I,S\}^\mathbb Z$.
\begin{thm}\label{thm:symdyn}
For every $(\alpha_i)_{i\in\mathbb Z}\in\Gamma$ there exists a
sequence $(x_i)_{i\in\mathbb Z}$ satisfying
\begin{enumerate}
\item $x_i\in |Q_{\alpha_i}|$ for $i\in\mathbb Z$,
\item $f_{\alpha_{i+1},\alpha_i}(x_i)=x_{i+1}$, for $i\in\mathbb Z$.
\end{enumerate}
Moreover, we have

\noindent{\bf periodic orbits:} if the sequence
$(\alpha_i)_{i\in\mathbb Z}$ is periodic with the principal period
$k$ then the trajectory $(x_i)_{i\in \mathbb Z}$ may be chosen so
that $x_k=x_0$, hence its trajectory is periodic

\noindent {\bf homo- and heteroclinic orbits:} if the sequence
$(\alpha_i)_{i\in\mathbb Z}$ is such that $\alpha_k=L_{i_-}$ for
$k\leq k_-$ and $\alpha_k=L_{i_+}$ for $k\geq k_+$, where
$i_-,i_+\in\{1,2\}$ then
\begin{equation*}
\lim_{k\to-\infty} x_k = L_{i_-}^*,\qquad \lim_{k\to\infty} x_k = L_{i_+}^*.
\end{equation*}
\end{thm}
\begin{proof}
The same as \cite[Theorem~7.1]{WZ}.
\end{proof}

\subsection{Resonance transitions.} Theorem~\ref{thm:symdyn}
implies the possibility for a comet to move between various
resonances. If we interpret staying close to $L_1^*$ or $L_2^*$
periodic orbits  as the $1:1$ resonance, then
Theorem~\ref{thm:symdyn} says that the comet can travel between
exterior and Sun regions in both directions and can move between
$5:3$, $3:2$, $1:2$, $2:3$ and $1:1$ resonances in an arbitrary
order.

\section{Symmetric periodic and homoclinic orbits.} \label{sec:symmetry}

In Section~\ref{subsec:symPCR3BP} the symmetry property of PCR3BP
and the associated Poincar\'e maps are described. In this section
we give the proof of the existence of an infinite number of
symmetric periodic and homoclinic orbits.
\begin{defn}
Let $I\ni t\to u(t)\in\Real^4$ be a solution of PCR3BP, where $I$
is the maximal interval of the existence of the solution. An orbit
$t\to u(t)$ is called \emph{$R$-symmetric}  iff
$$\mathrm{Image}(u)=\{u(t)\st t\in I\}=\{R(u(t))\st t\in
I\}=R(\mathrm{Image}(u)).$$
\end{defn}

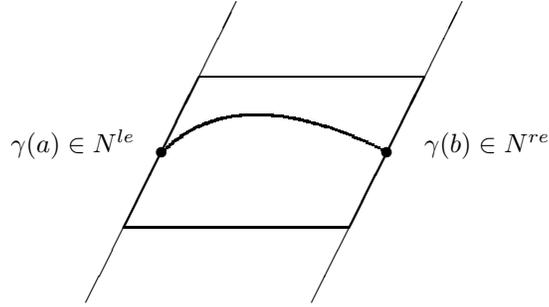
\begin{figure}[ptb]
\setlength{\unitlength}{1cm}
\begin{center}
\begin{picture}(7,4)
    \thinlines
    \put(4,0){\line(1,2){2}}
    \put(1,0){\line(1,2){2}}
    \put(1.5,1.){\line(1,0){3,0}}
    \put(2.5,3.){\line(1,0){3,0}}
    \thicklines
    \put(1.5,1.){\line(1,2){1}}
    \put(4.5,1.){\line(1,2){1}}
    \bezier{100}(2,2)(3,3)(5,2)
    \put(2,2){\circle*{.15}}
    \put(5,2){\circle*{.15}}
    \put(0,2){$\gamma(a)\in N^{le}$}
    \put(5.5,2){$\gamma(b)\in N^{re}$}
\end{picture}
 \caption{An h-set $N$ and a horizontal curve $\gamma$ in $N$\label{fig:h-curve}.}
\end{center}
\end{figure}

In this section we apply the method for finding symmetric
periodic, homo and heteroclinic orbits first introduced in
\cite{W2,W3} for the planar case  and later developed in
\cite{WZ1} in multidimensional situation. We recall here the basic
definitions.
\begin{defn}\label{def:h-curve}
Let $N$ be a h-set with one unstable and one stable direction and let
$\gamma:[a,b]\to\mathbb{R}^2$ be a continuous curve. We say that $\gamma$ is a horizontal
curve in $N$ if the following conditions hold:
\begin{enumerate}
 \item $\gamma((a,b))\subset \inter(|N|)$
 \item either $\gamma(a)\in N^{le} \quad and \quad \gamma(b)\in N^{re}$,\\
         or \quad $\mbox{ }\gamma(a)\in N^{re} \quad and \quad \gamma(b)\in N^{le}.$
\end{enumerate}
\end{defn}
The geometry of this concept is shown in Fig.~\ref{fig:h-curve}.
\begin{defn}\label{def:v-curve}
Let $N$ be a h-set with one unstable and one stable direction and let
$\gamma:[a,b]\to\mathbb{R}^2$ be a continuous curve. We say that $\gamma$ is a vertical
curve in $N$ if the following conditions hold:
\begin{enumerate}
 \item $\gamma((a,b))\subset \inter(|N|)$
 \item either $\gamma(a)\in N^{te} \quad and \quad \gamma(b)\in N^{be}$,\\
         or \quad $\mbox{ }\gamma(a)\in N^{te} \quad and \quad \gamma(b)\in N^{be}.$
\end{enumerate}
\end{defn}

The following theorem is a special case of \cite[Thm.3]{WZ1}.
\begin{thm}\label{thm:symmetry}
Assume $N_0,N_1,\ldots,N_k$ are h-sets with one unstable and one stable direction and
\begin{equation*}
N_0\bicover{f_0}N_1\bicover{f_1}\cdots\bicover{f_{k-1}}N_k.
\end{equation*}
If $\gamma:[a,b]\to\Real^2$ is a horizontal curve in $N_0$ and
$\bar\gamma:[\bar a,\bar b]\to\Real^2$ is a vertical curve in
$N_k$ then there exists $t_0\in(a,b)$ such that
\begin{eqnarray}
    (f_{m}\circ\cdots \circ f_0\circ\gamma)(t_0)\in \inter(|N_{m+1}|),\label{ass:1symmetry}
\end{eqnarray}
for $m=0,\ldots,k-1$ and
\begin{eqnarray}
    (f_{k-1}\circ\cdots \circ f_0\circ\gamma)(t_0)\in \bar\gamma((\bar a,\bar b)).\label{ass:2symmetry}
\end{eqnarray}
\end{thm}
Theorem~\ref{thm:symmetry} was first proven in \cite{W2} for a
planar case and direct covering relations. The generalization to a
higher dimension with one unstable direction and the direct
(forward) covering is presented in \cite{W3}. The proof of a
general situation (i.e. direct and backward covering in
multidimensional case) requires more sophisticated techniques and
is presented in \cite{WZ1}.

\subsection{Symmetric periodic orbits.}
In this section we will use Theorem~\ref{thm:symmetry} in order to
prove the existence of an infinite number of geometrically
different symmetric periodic orbits.

Before we state the main result in this section we introduce some
 notation. Let $(\alpha,\beta)\in\{L_1,L_2,X,E,I,S\}^2$ be an
admissible sequence of symbols. Let the maps $f_{\beta,\alpha}$ be
defined as in  (\ref{eq:def_f}).

\noindent{\bf Notation:} By
$Q_{\alpha}\bicover{f_{\beta,\alpha}}Q_{\beta}$ we will denote the
chain of covering relations associated with the sequence
$(\alpha,\beta)$, i.e.
$$Q_\alpha\bicover{P_1}V_1\bicover{P_2}V_2\bicover{P_3}\cdots
\bicover{P_{k-1}}V_{k-1} \bicover{P_k}Q_\beta,$$ where
$f_{\beta,\alpha}=P_k\circ\ldots\circ P_1$ and $V_i$,
$i=1,\ldots,k-1$ are suitable h-sets.

\begin{defn}
Let $f:X\to X$. By $\Fix(f)$ we will denote the set of fixed
points of $f$, i.e. $$\Fix(f)=\{y\in X\ |\ f(y)=y\}.$$
\end{defn}
\begin{figure}[htpb]
    \centerline{\includegraphics[width=2.5in]{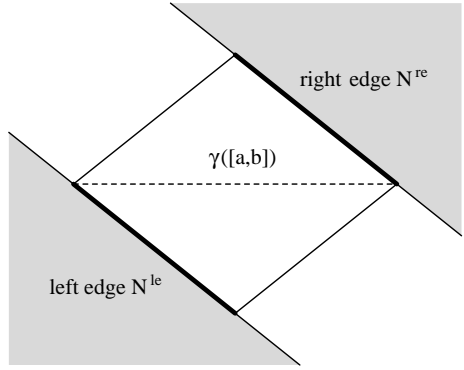}}
    \caption{A symmetric h-set. The $\gamma$ curve is both  horizontal
              and vertical  in $N$.\label{fig:symset}}
\end{figure}
\begin{thm}\label{thm:symperex} Let $\phi:\Real\times\Real^4\oarrow\Real^4$ denotes the
local flow induced by the PCR3BP with $C=3.03$ and
$\mu=0.0009537$. Assume $(\alpha_0,\alpha_1,\ldots,\alpha_n) \in
\{S,I,X,E,L_1,L_2\}^n$, $n>0$ is admissible sequence of symbols.
Then there exists a point $x_0\in |Q_{\alpha_0}|$ such that
\begin{eqnarray}
    (f_{\alpha_m,\alpha_{m-1}}\circ\cdots\circ f_{\alpha_1,\alpha_0})(x_0)\in |Q_{\alpha_m}|,\nonumber\\
    (f_{\alpha_{m-1},\alpha_m}^{-1}\circ\cdots\circ f_{\alpha_0,\alpha_1}^{-1})(x_0)\in
    |Q_{\alpha_m}|\label{eq:leftpart},
\end{eqnarray}
for $m=1,\ldots,n$, i.e., the trajectory of $x_0$ is coded by the
periodic sequence of symbols
\begin{equation}\label{eq:codeper}
(\alpha_n,\ldots,\alpha_1,\alpha_0,\alpha_1,\ldots,\alpha_n).
\end{equation}
Moreover,  $x_0$ is periodic and its orbit is $R$-symmetric.
\end{thm}
\begin{proof}
From the definitions of the h-sets used in the proof of homo- and
heteroclinic chains it follows that the sets
$H_1,H_2,V_0,G_0,E_0,F_0$ are $R$-symmetric. Therefore $\Fix(R)$
may be parameterized both as a horizontal and as a vertical curve
in each of these sets (see Fig~\ref{fig:symset}). Let
$\gamma:[a,b]\to |Q_{\alpha_0}|$ be the horizontal curve in
$Q_{\alpha_0}$ and $\bar\gamma:[\bar a, \bar b]\to|Q_{\alpha_n}|$
be the vertical curve in $Q_{\alpha_n}$, such that
$\gamma([a,b])\cup\bar\gamma([\bar a,\bar b])\subset \Fix(R)$.
Now, Theorem~\ref{thm:symmetry} applied to the sequence
\begin{equation*}
Q_{\alpha_0}\bicover{f_{\alpha_1,\alpha_0}}Q_{\alpha_1}\bicover{f_{\alpha_2,\alpha_1}}Q_{\alpha_2}\bicover{f_{\alpha_3,\alpha_2}}
\cdots\bicover{f_{\alpha_{n},\alpha_{n-1}}}Q_{\alpha_n}
\end{equation*}
implies that there exists a point $x_0=\gamma(t_0)\in
|Q_{\alpha_0}|\cap\Fix(R)$ such that
\begin{eqnarray*}
(f_{\alpha_m,\alpha_{m-1}}\circ\cdots\circ
f_{\alpha_1,\alpha_0})(x_0)\in |Q_{\alpha_m}|,\quad \text{for}
\quad m=1,\ldots,n,\\ (f_{\alpha_n,\alpha_{n-1}}\circ\cdots\circ
f_{\alpha_1,\alpha_0})(x_0) \in\bar\gamma((\bar a,\bar b))\subset
\Fix(R).
\end{eqnarray*}
From the definition of $f_{\beta,\alpha}$ (see
Eq.(\ref{eq:def_f})) as a composition of suitable Poincar\'e maps
it follows that there exists $T>0$ such that
\begin{equation*}
\phi(T,x_0)=(f_{\alpha_n,\alpha_{n-1}}\circ\cdots\circ
f_{\alpha_1,\alpha_0})(x_0)\in\Fix(R).
\end{equation*}
Since $R$ is the reversing symmetry of $\phi$ we obtain
\begin{equation*}
\phi(T,x_0)=R(\phi(T,x_0)) = \phi(-T,R(x_0)) = \phi(-T,x_0)
\end{equation*}
which proves  $x_0$ is periodic and its orbit  is  $R$-symmetric.

There remains to prove that the trajectory of $x_0$ is coded by the sequence
(\ref{eq:codeper}), i.e. (\ref{eq:leftpart}) is satisfied. We formulate this as a
separate lemma.
\end{proof}
\begin{lem}\label{lem:backtraj}
Assume $(\alpha_0,\ldots,\alpha_n)$ is an admissible sequence of
symbols. If $x\in\dom(f_{\alpha_n,\alpha_{n-1}}\circ\cdots\circ
f_{\alpha_1.\alpha_0})$ then
$R(x)\in\dom(f_{\alpha_{n-1},\alpha_n}^{-1}\circ\cdots\circ
f_{\alpha_0,\alpha_1}^{-1})$ and
\begin{equation*}
(R\circ f_{\alpha_m,\alpha_{m-1}}\circ\cdots\circ
f_{\alpha_1,\alpha_0})(x) =
(f_{\alpha_{m-1},\alpha_m}^{-1}\circ\cdots\circ
f_{\alpha_0,\alpha_1}^{-1}\circ R)(x),
\end{equation*}
for $m=1,\ldots,n$. Moreover, if $x=R(x)$ then
$$(f_{\alpha_{m-1},\alpha_m}^{-1}\circ\cdots\circ f_{\alpha_0,\alpha_1}^{-1})(x)\in |Q_{\alpha_m}|$$
for $m=1,\ldots,n$.
\end{lem}
\begin{proof}
One observes that if $(\alpha,\beta)$ is  admissible,  then
$(\beta,\alpha)$ is  admissible, too. Moreover, from the
definition of $f_{\beta,\alpha}$ (Eq. (\ref{eq:def_f})) it follows
that $R\circ f_{\beta,\alpha} = f^{-1}_{\alpha,\beta}\circ R$. Let
$x\in\dom(f_{\alpha_k,\alpha_{k-1}}\circ\cdots\circ
f_{\alpha_1.\alpha_0})$. Then
\begin{equation*}
\begin{split}
    (R\circ f_{\alpha_m,\alpha_{m-1}}\circ\cdots\circ f_{\alpha_1,\alpha_0})(x) = \\
    (f_{\alpha_{m-1},\alpha_m}^{-1}\circ R\circ
    f_{\alpha_{m-1},\alpha_{m-2}}\circ\cdots\circ f_{\alpha_1,\alpha_0})(x)=\\
    \cdots = (f_{\alpha_{m-1},\alpha_m}^{-1}\circ\cdots\circ f_{\alpha_0,\alpha_1}^{-1}\circ
    R)(x),
\end{split}
\end{equation*}
for $m=1,\ldots,n$. If in addition $x=R(x)$ then $x\in\dom(f_{\alpha_{m-1},\alpha_m}^{-1}\circ\cdots\circ
f_{\alpha_0,\alpha_1}^{-1})$ and
\begin{equation*}
(f_{\alpha_{m-1},\alpha_m}^{-1}\circ\cdots\circ
f_{\alpha_0,\alpha_1}^{-1})(x) \in
R(|Q_{\alpha_m}|)=|Q_{\alpha_m}|.
\end{equation*}
\end{proof}

\begin{rem}
Theorem~\ref{thm:symperex} implies that there exist infinitely
many geometrically different symmetric periodic orbits. This
follows immediately from the fact that there exists an infinite
number of admissible chains satisfying the assumptions of
Theorem~\ref{thm:symperex}.
\end{rem}

\subsection{Symmetric homoclinic orbits.}
In this section we apply Theorem~\ref{thm:symmetry} in order to
prove the existence of  infinitely many geometrically different
symmetric homoclinic orbits to $L_1^*$ and $L_2^*$ Lyapunov
orbits.

The following theorem shows how to use the method of covering
relations in order to prove the existence of symmetric homoclinic
or heteroclinic orbits. Later we will apply it to Poincar\'e maps
for PCR3BP.
\begin{thm}\label{thm:symhom}
Let $N_0,N_1,\ldots,N_k$ be h-sets, such that
\begin{equation*}
N_0\bicover{f_0}N_1\bicover{f_1}\cdots\bicover{f_{k-1}}N_k\bicover{f_k}N_k
\end{equation*}
and let $\gamma:[a,b]\to\Real^2$ in $N_0$ be a horizontal curve in
$N_0$. If $f_k$ is hyperbolic (see \cite[Def. 1]{gaz}) on $N_k$,
then there exists a point $x_0\in \gamma((a,b))$ such that
\begin{eqnarray*}
    (f_m\circ \cdots \circ f_0)(x_0) \in \inter(|N_{m+1}|),\quad \text{for }  m=0,1,\ldots,k-1,\\
    (f_k^n\circ f_{k-1}\circ\cdots\circ f_0)(x_0)\in \inter (|N_k|),\quad \text{for } n>0.
\end{eqnarray*}
Moreover,
\begin{eqnarray*}
    \lim_{n\to\infty}(f_k^n\circ f_{k-1}\circ\cdots\circ f_0)(x_0)=x_*,
\end{eqnarray*}
where $x_*$ is a unique fixed point of $f_k$ in $|N_k|$.
\end{thm}
\begin{proof}
From Theorem~\ref{thm:symmetry} it follows that for every $n>0$
there exists $t_n\in[a,b]$ such that
\begin{eqnarray*}
    (f_m\circ \cdots \circ f_0)(\gamma(t_n)) \in \inter(|N_{m+1}|),\quad \text{for } m=0,1,\ldots,k-1,\\
    (f_k^n\circ f_{k-1}\circ\cdots\circ f_0)(\gamma(t_n))\in \inter (|N_k|).
\end{eqnarray*}
Since $\gamma([a,b])$ is compact we can find $t_*\in[a,b]$ such that
\begin{eqnarray*}
    (f_m\circ \cdots \circ f_0)(\gamma(t_*)) \in \inter(|N_{m+1}|),\quad \text{for } m=0,1,\ldots,k-1,\\
    (f_k^n\circ f_{k-1}\circ\cdots\circ f_0)(\gamma(t_*))\in \inter (|N_k|),\quad \text{for } n>0.
\end{eqnarray*}
Since neither $f(\gamma(a))\notin N_1$ nor $f(\gamma(b))\notin
N_1$ we get $t_*\in(a,b)$. Now, $f_k$ is hyperbolic on $N_k$.
Therefore by Theorem 3 in \cite{gaz},
\begin{equation*}
    \lim_{n\to\infty}(f_k^n\circ f_{k-1}\circ\cdots\circ f_0)(x_0)=x_*,
\end{equation*}
where $x_0:=\gamma(t_*)$.
\end{proof}

Now we can state the basic result in this section.
\begin{thm}\label{thm:symhomex}
Assume $(\alpha_0,\alpha_1,\ldots,\alpha_n)$ is an admissible
nonconstant chain of symbols $\{S, I, X, E, L_1, L_2 \}$, such
that $\alpha_n\in\{L_1,L_2\}$. Then there exists a symmetric
homoclinic orbit associated with the sequence of symbols
\begin{equation}\label{eq:codehom}
(\ldots,\alpha_n,\alpha_n,\alpha_{n-1},\ldots,\alpha_1,\alpha_0,\alpha_1,\ldots,\alpha_{n-1},\alpha_n,\alpha_n,\ldots).
\end{equation}
\end{thm}
\begin{proof}
Let $\gamma:[a,b]\to|Q_{\alpha_0}|$ be a horizontal curve in
$Q_{\alpha_0}$ such that $\gamma([a,b])\subset\Fix(R)$. From
Lemma~5.5 in \cite{WZ} it follows that $P_+$ is hyperbolic on
$|H_1|=Q_{L_1}$ and $P_-$ is hyperbolic on $|H_2|=Q_{L_2}$. Since
$\alpha_n\in\{L_1,L_2\}$ Theorem~\ref{thm:symhom} there exists
$x_0\in\gamma((a,b))$ such that
\begin{eqnarray*}
    (f_{\alpha_m,\alpha_{m-1}}\circ \cdots \circ f_{\alpha_1,\alpha_0})(x_0) \in \inter(|Q_{\alpha_m}|),\quad \text{for } m=1,\ldots,n,\\
    (f^k_{\alpha_n,\alpha_n}\circ f_{\alpha_n,\alpha_{n-1}}\circ \cdots \circ f_{\alpha_1,\alpha_0})(x_0)\in \inter (|Q_{\alpha_n}|),\quad \text{for }
    k>0,\\
    \lim_{k\to\infty}(f^k_{\alpha_n,\alpha_n}\circ f_{\alpha_n,\alpha_{n-1}}\circ \cdots \circ
    f_{\alpha_1,\alpha_0})(x_0) = L,
\end{eqnarray*}
where $L=L_1^*$ or $L=L_2^*$ is a unique fixed point in $|Q_{\alpha_n}|$. Since $x_0=R(x_0)$ Lemma~\ref{lem:backtraj} implies that
\begin{eqnarray*}
    (f_{\alpha_{m-1},\alpha_{m}}^{-1}\circ \cdots \circ f_{\alpha_0,\alpha_1}^{-1})(x_0) \in \inter(|Q_{\alpha_m}|),\quad \text{for } m=1,\ldots,n,\\
    (f^{-k}_{\alpha_n,\alpha_n}\circ f_{\alpha_{n-1},\alpha_{n}}^{-1}\circ \cdots \circ f_{\alpha_0,\alpha_1}^{-1})(x_0)\in \inter (|Q_{\alpha_n}|),\quad \text{for }
    k>0,\\
    \lim_{k\to\infty}(f^{-k}_{\alpha_n,\alpha_n}\circ f_{\alpha_{n-1},\alpha_{n}}^{-1}\circ \cdots \circ f_{\alpha_0,\alpha_1}^{-1})(x_0) =
    R(L)=L
\end{eqnarray*}
which proves that the trajectory of $x_0$ is a symmetric homoclinic orbit coded by the
sequence of symbols (\ref{eq:codehom}).
\end{proof}

\begin{rem}
Theorem~\ref{thm:symhomex} implies that there exist infinitely
many symmetric homoclinic orbits which are geometrically
different. This follows immediately from the fact that there
exists an infinite number of admissible chains satisfying the
assumptions of Theorem~\ref{thm:symhomex}.
\end{rem}

\section{Technical data.}
The computer assisted proofs of Lemma~\ref{lem:cov_ext_reg} and
Lemma~\ref{lem:cov_int_reg} will be not discussed here. All ideas
involved in such proof  were presented in \cite{WZ}. The C++
sources containing the rigorous numerical proof of  Lemmas from
\cite{WZ}, Lemma~\ref{lem:cov_ext_reg} and
Lemma~\ref{lem:cov_int_reg} is available at \cite{W1}.

The program uses the interval arithmetic and set algebra package
developed at Jagiellonian University by CAPD group \cite{CAPD}.

\section{Appendix. Reading resonances from the trajectory in rotating frame.}
We describe the heuristic approach, which allows to read the
resonance from the inspection of the trajectory in the rotating
coordinate frame.

We assume that Jupiter and comet move in the heliocentric inertial
frame in the counterclockwise direction and the distance comet-Sun
has well visible maxima or minima along the trajectory. This means
that an approximate ellipse on which the comet is moving has
nonzero eccentricity.

Let $R$ denote the resonance. Let $T$ be an approximate period of
the comet in the heliocentric frame. Let us recall that the period
of the Jupiter is equal to $1$. Hence
\begin{equation}
   R = \frac{1}{T}
\end{equation}
Then in the heliocentric inertial frame the average angular
velocity of the comet is $\frac{2 \pi}{T}$ and that of the Jupiter
is equal to $2 \pi$.

For an approximate periodic trajectory of a comet in the rotating
frame we introduce the following notation
\begin{itemize}
\item $\theta$  is the number of full turns around the Sun during
the whole period. This number is positive for trajectories in the
interior region and negative in the exterior region.
\item $M$ - the number of maxima (or minima) of the distance
between the Sun and the comet.
\end{itemize}
Since the distance Sun-comet reaches the maximum (or minimum) only
when the comet is at the aphelion (or perihelion), hence
consecutive maxima (minima) occur with the period $T$. In the
rotating frame the difference between the angular variables of the
comet and Jupiter is equal to $2\pi \theta/ M$. Observe that this
difference is the same in  both reference frames, the  inertial
one and the rotating one. Hence
\begin{eqnarray*}
   \frac{2 \pi \theta}{M} &=& \left(\frac{2 \pi}{T} - 2 \pi \right)
   T \\
   \frac{\theta}{M}&=& 1 - T \\
   T&=& \frac{M - \theta}{M}.
\end{eqnarray*}
Hence finally
\begin{equation}
   R = \frac{M}{M - \theta}.  \label{eq:resonance}
\end{equation}
Let us apply (\ref{eq:resonance}) to Figures~\ref{fig:homoc1} and
\ref{fig:homoc2}. For interior homoclinics we count the maxima and
for exterior homoclinic we count the minima. We have
\begin{itemize}
\item  the interior homoclinic orbit in Figure~\ref{fig:homoc1}:
$\theta=1$, $M=3$. Hence $R = \frac{3}{2}=3:2$.
\item  the interior homoclinic orbit in Figure~\ref{fig:homoc2}:
$\theta=2$, $M=5$. Hence $R = \frac{5}{3}=5:3$.
\item  the exterior homoclinic orbit in Figure~\ref{fig:homoc1}:
$\theta=-1$, $M=1$. Hence $R = \frac{1}{2}=1:2$.
\item  the exterior homoclinic orbit in Figure~\ref{fig:homoc2}:
$\theta=-1$, $M=2$. Hence $R = \frac{2}{3}=2:3$.
\end{itemize}

\end{document}